\documentclass[12pt,letterpaper]{article}
\usepackage[centertags]{amsmath}
\usepackage{amsfonts}
\usepackage{amssymb}
\usepackage{mparhack,makeidx}
\usepackage{amsthm}
\usepackage{color}
\usepackage[T1]{fontenc}
\usepackage[leftcaption]{sidecap}
\usepackage{subfigure}
\usepackage[centertags]{amsmath}
\usepackage{amssymb}
\usepackage{amsthm}
\usepackage[active]{srcltx}
\usepackage[latin1]{inputenc}
\usepackage[dvips]{graphicx}
\usepackage{amsthm,amsfonts} 
\usepackage{graphicx} 

\newtheorem{thm}{Theorem}[section]
\newtheorem{cor}[thm]{Corollary}
\newtheorem{lem}[thm]{Lemma}
\newtheorem{prop}[thm]{Proposition}
\newtheorem{defn}[thm]{Definition}
\newtheorem{exe}[thm]{Example}

\newtheorem{remark}[thm]{Remark}
\newtheorem{conjecture}[thm]{Conjecture}
\newcommand{\Ge}{\|g\|^b}
\newcommand{\Real}{\mathbb{R}}
\newcommand{\Complex}{\mathbb{C}}

\newcommand{\tl}{\tilde}

\newcommand{\Nat}{\mathbb{N}}

\newcommand{\Int}{\mathbb{Z}}

\newcommand{\hs}{\hspace{.5cm}}
\newcommand{\vs}{\vspace{.5cm}}
\newcommand{\nd}{\noindent}
\newcommand{\ds}{\displaystyle\frac}
\newcommand{\dss}{\displaystyle}

\newcommand{\hp}{$H^p(\mathbb{D}, F)$}
\newcommand{\bl}{\blacksquare}

\newcommand{\al}{\alpha}

\newcommand{\p}{{\nd \bf{Proof:\hspace{.2cm}}}}
\newcommand{\holo}{{\cal{H}}(\mathbb{D};F)}
\newcommand{\hol}{{\cal{H}}(\mathbb{D};F')}

\newcommand{\hf}{$H^{\infty}(\mathbb{D}, F)$}

\makeindex

\begin{document}

\title{Duality of certain Banach spaces of vector-valued
holomorphic functions}
\author{F. J. Bertoloto\thanks{Supported by FAPESP.\hfill\newline Mathematics Subject Classification: 46G20, 46G10, 30H10.\hfill\newline Keywords: Hardy spaces; holomorphic functions; RNP, ARNP and UMDP properties; duality results; weak version of Banach spaces.}}
\date{}
\maketitle

\nd Faculdade de Matemática\\
\nd Universidade Federal de Uberlândia\\
\nd 38.400-902 - Uberlândia - Brazil\\
\nd e-mail: bertoloto@famat.ufu.br

\section*{Abstract}

\numberwithin{equation}{section}

\hs In this work we study the vector-valued Hardy spaces
$H^p(\mathbb{D};F)$ ($1\leq p \leq \infty$) and their relationship with RNP, ARNP and the UMDP properties.
By following the approach of Taylor \cite{Taylor 1},\cite{Taylor
2} in the scalar-valued case, we prove that, when $F$ and $F'$
have the ARNP property, then $H^p(\mathbb{D};F)'$ and
$H^q(\mathbb{D};F')$ are canonically topologically isomorphic (for
$p,q\in (1,\infty)$ conjugate indices) if and only if $F$ has the UMDP.

\section{Introduction}

Throughout this paper let $\mathbb{D}=\{z\in \Complex; |z|<1\}$, $\mathbb{T}=\{z\in\Complex
; |z|=1\}$, and let $F$ be a complex Banach space. Let $\mathcal{H}(\mathbb{D};F)$ denote the vector space of all holomorphic functions $f\colon\mathbb{D}
\longrightarrow F$.
Let $1\leq p,q\leq \infty$ with $\frac{1}{p}+\frac{1}{q}=1$.

This paper is devoted to the study of subspaces of $\mathcal{H}(\mathbb{D};F)$ that satisfy certain axioms. Our paper is indeed a vector-valued versions of
Taylor \cite{Taylor 1, Taylor 2}. Our main objective is the study of the vector-valued Hardy space $H^p(\mathbb{D};F)$. In particular we show that if $F$ and $F^\prime$ have
the analytic Radon-Nikodym property (ARNP) and $1<p<\infty$, then the spaces $H^p(\mathbb{D};F)^\prime$ and  $H^q(\mathbb{D};F^\prime)$ are canonically topologically isomorphic
if and only if $F$ has the unconditional martingale difference property (UMDP).

This paper is organized as follows. Section 2 is devoted to the study of subspaces of  $\mathcal{H}(\mathbb{D};F)$ satisfying certain axioms. In particular we study the spaces
 $H^p(\mathbb{D};F)$. Section 3 is devoted to the study of the space  $H^p(\mathbb{T};F)$ of boundary-value functions associated with  $H^p(\mathbb{D};F)$. There we introduce
the ARNP following Bukhvalov \cite{Bukhvalov}. Section 4 is devoted to the study of duality theory for the spaces  $H^p(\mathbb{D};F)$, and its connection with the UMDP. Finally
Section 5 is devoted to the study of the spaces  $H^p(\mathbb{D};F)_w$ and  $L^p(\mathbb{T};F)_w$, the weak versions of the spaces  $H^p(\mathbb{D};F)$ and  $L^p(\mathbb{T};F)$.

This paper is based on the doctoral thesis of the author at UNICAMP, written under the supervision of Professor Jorge Mujica.

\section{Preliminaries}

For a few elementary facts regarding the Bochner integral we refer to Mujica \cite{Mujica}.
Except for Proposition \ref{teo12}, that improves a result of Taylor \cite[Theorem 8.1]{Taylor 1}, all the other results have proofs
similar to those given by Taylor \cite{Taylor 1, Taylor 2} in the scalar-valued case. We cite
Lemma \ref{nonincreasinglemma} with a proof in order to highlight its importance for other results.
Here, we change the notation of Taylor trying to adapt to a more current one.

\subsection{Vector-valued holomorphic functions on the open unit disc}

\hs It is well known that if $f\in \mathcal{H}(\mathbb{D};F)$ then

$$f(z)=\sum_{n=0}^{\infty}\gamma_n(f)z^n,$$

\nd where $\gamma_n(f)\in F$. This series converges uniformly and
absolutely in the open disks $D_r=\{ z;|z|<r\}$ for $0< r<1$. It is
easy to see that

$$\gamma_n(\al f+\beta g)=\al\gamma_n(f)+\beta\gamma_n(f),$$

\nd for every $f,g\in {\cal{H}}(\mathbb{D};F)$ and $\alpha, \beta\in
\Complex$.

\begin{defn}{\rm Let $z\in\mathbb{D}$.

\nd{\bf a)} If $v\in F$ and $\|v\|=1$, then we define
$u^v_n(z)=z^nv$.

\nd{\bf b)} \label{def 15}If $f\in {\cal{H}}(\mathbb{D};F)$ and $t\in
\Real$, then we define $U_tf(z)=f(ze^{it})$.

\nd{\bf c)} \label{def 16}If $f\in {\cal{H}}(\mathbb{D};F)$ and $|w|\leq
1$, then we write $T_wf(z)=f(zw)$.}
\end{defn}

\hs The reader can verify that $U_t=T_{e^{it}}$ and
$\gamma_n(T_wf)=w^n\gamma_n(f)$ for every $t\in \Real$, $w\in
\mathbb{D}$ and $n\in \Nat$. Besides, $U_t$ and $T_w$ are linear
operators on ${\cal{H}}(\mathbb{D};F)$.

\begin{defn}\label{def2}
{\rm Let $f\in {\cal{H}}(\mathbb{D};F)$ and $g\in
{\cal{H}}(\mathbb{D};F')$. Then:

\begin{center}
$f(z)=\dss\sum_{n=0}^{\infty}\gamma_n(f)z^n$ \hs and\hs
$g(z)=\dss\sum_{n=0}^{\infty}\gamma_n(g)z^n$.
\end{center}

\hs For each $z\in \mathbb{D}$, let

$$B(f,g;z)=\dss\sum_{n=0}^{\infty}<\gamma_n(f),\gamma_n(g)>z^n,$$

\nd where $<\gamma_n(f),\gamma_n(g)>=\gamma_n(g)(\gamma_n(f))$.
}
\end{defn}

\begin{prop}\label{propp1}{\rm Let $f\in\mathcal{H}(\mathbb{D};F)$, $g\in\mathcal{H}(\mathbb{D};F')$, $w\in\mathbb{D}$ and $z\in \mathbb{D}$. Then,

\nd{\bf a)} $B(f,g;.)\in{\cal{H}}(\mathbb{D})$.

\nd{\bf b)} $B(f,g;z)$ is linear on $f$ for each $g$ and $z$.

$B(f,g;z)$ is linear on $g$ for each $f$ and $z$.

\nd{\bf c)} $B(T_wf,g;z)=B(f,g,wz)$.

\nd{\bf d)}
$B(f,g;z)=\ds{1}{2\pi}\int_0^{2\pi}<f(z_1e^{i\theta}),g(z_2e^{-i\theta})>d\theta$
whenever $z=z_1z_2$, with $z_1,z_2\in \mathbb{D}$.
}
\end{prop}

\hs Let $H\subset{\cal{H}}(\mathbb{D};F)$ be a complex normed
vector space. We say $H$ {\it is a normed space of type
${\cal{H}}(\mathbb{D};F)$}, or is {\it of type
${\cal{H}}(\mathbb{D};F)$}, if it contains at least two elements. If
$H$ is a Banach space, we say $H$ is {\it a Banach space of type
${\cal{H}}(\mathbb{D};F)$}.

\hs Let $A>0$ and $n\in \Nat$. We introduce now seven properties
(or axioms) that a space $H$ of type ${\cal{H}}(\mathbb{D};F)$ can
satisfy.

\vs

\nd {$P_1$:} There exists a constant $A$ such that
$\|\gamma_n(f)\|\leq A\|f\|$ if $f\in H$. The least such $A$ we denote by $A_1(H)$. Thus each $\gamma_n\colon H\longrightarrow F$
is a continuous linear operator and $\|\gamma_n\|\leq A_1(H)$.

\vs \nd{\bf $P_2$:} $u_n^v\in H$ and there exists a constant $A$
such that $\|u_n^v\|\leq A$ for all $v\in F$ with $\|v\|=1$. The
least such $A$ we denote by $A_2(H)$. Thus
$\|u_n^v\|\leq A_2(H)$.

\vs

\hs For the axioms $P_3$ and $P_4$, we need Definition
\ref{def 15}.

\vs

\nd{\bf $P_3$:} $U_tf\in H$ if $f\in H$ and $t\in \Real$. Besides
$\|U_tf\|=\|f\|$ and $U_t$ is an isometry.

\vs

\nd{\bf $P_4$:} If $f\in H$ and $0<r<1$, then $T_rf\in H$. There
exists a constant $A$ such that $\|T_rf\|\leq A\|f\|$. The least
such $A$ we denote by $A_4(H)$. Thus $\|T_r\|\leq A_4(H)$.

\vs

\nd $P_5$: If $f\in H$ and $0<r<1$, then $T_rf\in H$ and
$\|f\|=\dss\sup_{0\leq r<1}\|T_rf\|$.

\vs

\nd $P_6$: If $f\in H$ and $0<r<1$, then $T_rf\in H$ and
$\dss\lim_{r\rightarrow 1}\|T_rf-f\|=0$.

\vs

\nd $P_7$: If $f\in \holo$ is such that $T_rf\in H$ whenever
$0<r<1$ and $\dss\sup_{0\leq r<1}\|T_rf\|<\infty$, then $f\in H$ and
$\|f\|=\displaystyle\sup_{0\leq r<1}\|T_rf\|$.

\hs We write $H'$ for the topological dual of $H$. We can express
the constants $A_1(H),A_2(H)$ and $A_4(H)$ as follows:

\begin{center}
$A_1(H)=\dss\sup_{n\in\Nat}\|\gamma_n\|$,\hs
$A_2(H)=\hspace{-.6cm}\dss\sup_{\begin{array}{cc}
\vspace{-.7cm}\\_{v\in F,\|v\|=1}\vspace{-.2cm}\\
_{n\in\Nat}\end{array}}\hspace{-.6cm}\|u_n^v\|$, \hs
$A_4(H)=\dss\sup_{0\leq r<1}\|T_r\|$.

\end{center}

\begin{prop}\label{lem41} {\rm Let $H$ be of type
$\mathcal{H}(\mathbb{D};F)$. Then

\nd {\bf a)} If $H$ satisfies $P_1$ e $P_2$, then $A_1(H)A_2(H)\geq 1.$

\nd {\bf b)} If $H$ satisfies $P_4$ and for some $n\in \Nat$,
$u_n\in H$ (in particular, if $H$ satisfies $P_2$), then
$A_4(H)\geq 1$.

\nd {\bf c)} $P_5$ implies $P_4$ and $A_4(H)=1$.

\nd {\bf d)} $P_4$ and $P_7$ imply $P_5$.}
\end{prop}
\p For assertion (a) we observe that $\gamma_n(u_n^v)=v$. Then if
the properties {\bf $P_1$} and {\bf $P_2$} are satisfied, it
follows that:

$$1=\|\gamma_n(u_n^v)\|\leq A_1(H)A_2(H).$$

\hs For assertion (b), we observe that
$T_r(u_n^v)=r^nu_n^v$. Assertions (c) and (d) are not difficult to prove. $\bl$

\hs A given normed space $H$ is said to be a {\it  normed space of
type ${\cal{H}}(\mathbb{D};F)_k$}, if it is a normed space of type
${\cal{H}}(\mathbb{D};F)$ that satisfies $P_1,\ldots, P_k$
($k=1,\ldots,7$). We say {\it Banach space of type
${\mathcal{H}}(\mathbb{D};F)_k$} if $H$ is a Banach space.

\begin{defn}{\rm Let $H^p(\mathbb{D};F)$ denote the space of all $f\in \mathcal{H}(\mathbb{D};F)$
such that $\|f\|_p< \infty$, where

\begin{center} $\|f\|_p=\dss\sup_{0\leq
r<1}\left(\ds{1}{2\pi}\int_0^{2\pi}
\|f(re^{i\theta})\|^pd\theta\right)^{\frac{1}{p}}$ if $1\leq
p<\infty$
\end{center}

and

$$\|f\|_\infty=\dss\sup_{z\in \mathbb{D}}\|f(z)\|.$$

 Then $H^p(\mathbb{D};F)$ is a Banach space for the norm $\|.\|_p$. The proof
 of Rudin
\cite[page 331]{Rudin} in the scalar-valued case applies. When
$F=\Complex$, we write $H^p(\mathbb{D})$ instead of
$H^p(\mathbb{D};\Complex)$.}
\end{defn}

\hs Next we will see that the spaces $H^p(\mathbb{D};F)$ are
examples of Banach spaces of type $\mathcal{H}(\mathbb{D};F)_4$.
To see this we need the following lemma.

\begin{lem}\label{nonincreasinglemma}{\rm For each $0\leq r<1$, $1\leq p\leq \infty$ and $f\in
{\mathcal{H}}(\mathbb{D};F)$, let

$$\|f\|_{p,r}=\left(\ds{1}{2\pi}\int_0^{2\pi}
\|f(re^{i\theta})\|^pd\theta\right)^{\frac{1}{p}}
$$

and

$$\|f\|_{\infty,r}=\dss\max_{0\leq \theta\leq 2\pi}\|f(re^{i\theta})\|.$$

Then the maps $r\mapsto \|f\|_{p,r}$ and $p\mapsto \|f\|_{p,r}$
are nondecreasing functions of $r$ and $p$ respectively. Also,

\begin{equation}\label{sup}\dss\lim_{p\rightarrow
\infty}\|f\|_{p,r}=\|f\|_{\infty,r}, \hs \dss\lim_{r\rightarrow
1}\|f\|_{p,r}=\|f\|_{p} .\end{equation}}\end{lem}

\p Let $1\leq p<p_1$. We consider the functions

\begin{equation}\label{pincreasing}\phi=\|f\|^p\in L^{\frac{p_1}{p}}(T)\hs {\rm and}\hs
\psi=1\in L^{\frac{p_1}{p_1-p}}.
\end{equation}

Now is enough to apply the Holder's inequality.

The Holder's inequality also resolves the case for $r$, when
$1\leq p<\infty$, with the help of the Cauchy's Integral Formula
(see, e.g., \cite[page 8]{Duren}). When $p=\infty$ we use the
Maximum Modulus Theorem (see, e.g., Mujica \cite{Mujica} or Thorp
\cite[page 641]{Thorp})).

Concerning the first equality in (\ref{sup}), we observe that the
limit exists, since $p\mapsto \|f\|_{p,r}$ is a
non-decreasing function of $p$ that is bounded by
$\|f\|_{\infty,r}$.

Given $\epsilon>0$, we can assure there exists a set
$A\subset[0,2\pi]$ of positive Lebesgue measure $\rho\leq 1$
satisfying,

$$\|f(re^{i\theta})\|>\|f\|_{\infty,r}-\epsilon$$

\nd for every $\theta\in A$. Hence,

$$\|f\|_{p,r}\geq \left(\ds{1}{2\pi}\int_A\|f(re^{i\theta})\|^pd\theta\right)^{\frac{1}{p}}\geq
(\|f\|_{\infty,r}-\epsilon)(\rho)^{\frac{1}{p}}.$$

Taking $p\longrightarrow \infty$, the lemma is proved, since the second
equality in (\ref{sup}) is clear.$\bl$

\begin{prop}\label{prop31}{\rm If $1\leq p\leq\infty$, then \hp\, is a
Banach space of type $\holo_4$. Also, $A_k($\hp$)=1$ (
$k=1,2,4$)}.
\end{prop}

\p By Lemma \ref{nonincreasinglemma} and the proof in the scalar-valued case, the proposition follows.
See Taylor \cite[Theorem 11.1]{Taylor 2}.$\bl$

\begin{remark}{\rm The following equalities hold true for $f\in
\holo$ and $\rho>0$ such that $0\leq r<\rho<1$:

\begin{equation}\label{ig1} \|T_rf\|_{p,\rho}=\left(\ds{1}{2\pi}\int_{0}^{2\pi}\|T_rf(\rho
e^{i\theta})\|^pd\theta\right)^{\frac{1}{p}}=\left(\ds{1}{2\pi}\int_{0}^{2\pi}\|f(r\rho
e^{i\theta})\|^pd\theta\right)^{\frac{1}{p}}=\|f\|_{p,r\rho}.
\end{equation}

Thus, from equation (\ref{ig1}), for each
$f\in\holo$ we have

\begin{equation}\label{ig2}
\|T_rf\|_p=\|f\|_{p,r}.
\end{equation}}
\end{remark}

 This happens since $r\longmapsto \|f\|_{p,r}$ is a
nondecreasing function of $r$ (Lemma \ref{nonincreasinglemma}).

\begin{prop}\label{ig04}{\rm The space \hf\, is a Banach space of type $\mathcal{H}(\mathbb{D};F)_5$. It also satisfies $P_7$
but not $P_6$.}
\end{prop}

\p The proof in the scalar-valued case applies. See Taylor \cite[Theorem 12.1]{Taylor 2}. $\bl$

\subsection{Additional Properties of the normed spaces of
type $\mathcal{H}(\mathbb{D};F)_k$}

\begin{prop}\label{teo1}{\rm If $H$ is a normed space of type $\holo_1$, $f\in
H$ and $z\in \mathbb{D}$, then:

$$\|f(z)\|\leq \ds{A_1(H)\|f\|}{1-|z|}.$$}

\end{prop}

\p By the expansion $f(z)=\sum_{n=0}^{\infty}\gamma_n(f)z^n,$ the result follows. $\bl$

\begin{defn}\label{defff1}{\rm Let $H$ be of type $\holo_1$, $f\in H$, $g\in {\cal{H}}(\mathbb{D}, F')$ and $z\in
\mathbb{D}$. Then

$$N(g;z)=\dss\sup_{\|f\|=1}|B(f,g;z)|.$$}

\end{defn}

\nd $N(g;z)$ has the following properties:

\nd{\bf a)} $N(g+h;z)\leq N(g;z)+N(h;z)$;

\nd{\bf b)} $N(ag;z)=|a|N(g;z)$;

\nd{\bf c)} $N(T_wg;z)=N(g;wz)$, $|w|\leq 1$.

\nd{\bf d)} $N(g;z)=N(g;|z|)$.

\nd{\bf e)} If $N(g;r)=N(g;z)$, where $|z|=r$, then $N(g;r)$ is a continuous and increasing function of $r$.

\nd{\bf f)} $N(g;r)=0$, for $0\leq r<1$,

\hs The proof of properties (a), (b) and (c) is easy. We omit it. For item (d), (e) and (f), see Taylor \cite[Theorem 6.1]{Taylor 1}.

\begin{prop}\label{teo 3} {\rm Let $H$ be a Banach space of type
$\holo_2$, $f\in \holo$ and $w\in \mathbb{D}$. Then $T_wf\in H$
and the function $w\longmapsto T_wf$ is analytic on $\mathbb{D}$
with series expansion:

\begin{equation}\label{eq1}
T_wf=\dss\sum_{n=0}^{\infty}a_n.u_n^{v_n} ,
\end{equation}

\nd where

\nd{\bf a)} $a_n=\|\gamma_n(f)\|w^n$ e
$v_n=\ds{\gamma_n(f)}{\|\gamma_n(f)\|}$, if $\gamma_n(f)\neq 0$;

\nd{\bf b)} $a_n=0\in \Complex$ e $v_n=0\in F$, if $\gamma_n(f)= 0$. }
\end{prop}

\p The proof in the scalar-valued case applies. See Taylor \cite[Theorem 4.1]{Taylor 1}. $\bl$

\subsection{Holomorphic functions with values in a dual Banach space}

 \hs In this section we will always consider spaces $H$ of type
$\holo_3$, and from them, we will define a space of type $\hol_4$
that we call $H^b$.

\begin{defn}{\rm Let $H$ be of type $\holo_3$. Define $H^b$ as the set of $g\in {\cal H}(\mathbb{D},F')$ such that $N(g;r)$ is a
bounded function of $r$. We write for $g\in H^b$,

\begin{center}$N(g)=\dss\lim_{r\rightarrow 1}N(g;r)$\end{center}}

\end{defn}

\nd where the limit can be replaced by $N(g)=\dss\sup_{0\leq
r<1}N(g;r)$ from property (e) after Definition \ref{defff1}.

\begin{prop}{\rm If $H$ is of type $\holo_3$, then
$H^b$ is a Banach space with norm $N(g)$.}
\end{prop}

\p The proof in the scalar-valued case applies. See Taylor \cite[Theorem 7.1]{Taylor 1}. $\bl$

\vs

\nd {\bf Notation:} From now on $N(g)$ will be denoted by $\Ge$, for
all $g\in H^b$.

\begin{prop}\label{prop4.3}{\rm Let $H$ be of type $\holo_3$. Then
$H^b$ is a Banach space of type $\hol_4$. The constants
$A_k(H^b)$, $k=1,2,4$, satisfy

\nd{\bf a)} $A_1(H^b)\leq A_2(H)$.

\nd{\bf b)} $A_2(H^b)\leq A_1(H)$.

\nd{\bf c)} $A_4(H^b)=1$.
}
\end{prop}

\p See the scalar-valued case in Taylor \cite[Theorem 7.2]{Taylor 1}. The proof here is similar. $\bl$

\begin{prop} \label{teo12}{\rm Let $H$ be a space of type $\mathcal{H}(\mathbb{D};F)_4$.
Let $\gamma\in H'$ and let $g\colon \mathbb{D}\longrightarrow F'$ be defined by

\begin{equation}
\label{eq10}g(z)=\dss\sum_{n=0}^{\infty}z^n{\gamma_n}
\end{equation}

\nd for all $z\in \mathbb{D}$, where $\gamma_n(v):=\gamma(u_n^v)$ for
every $v\in F$. Then $g\in H^b$ and $\Ge\leq A_4(H)\|\gamma\|$.}

\end{prop}

\p Firstly, we have $v \neq 0$

$$\begin{array}{rl}|\gamma_n(v)|=|\gamma(u_n^v)| &\leq
\|\gamma\|\|u_n^v\|\vspace{.03cm}\\ &
=\|\gamma\|.\|v\|.\|u_n^{\frac{v}{\|v\|}}\| \vspace{.03cm}\\
& \leq \|\gamma\|A_2(H)\|v\|.
\end{array}$$

\hs Thus, $\|\gamma_n\|\leq \|\gamma\|A_2(H)$ which implies
$\gamma_n\in F'$ and $g\in\hol$. By Proposition \ref{teo
3}, we have for $w\in \mathbb{D}$ and $f\in H$ that

$$T_wf=\dss\sum_{n=0}^{\infty}a_n.u_n^{v_n},$$

\nd is a element of $H$, where

\nd i) $a_n=\|\gamma_n(f)\|w^n$ and
$v_n=\ds{\gamma_n(f)}{\|\gamma_n(f)\|}$, if $\gamma_n(f)\neq 0$.

\nd ii) $a_n=0\in \Complex$, if $\gamma_n(f)= 0$.

\hs Let $\gamma\in H'$ and $\{a_{n_k}\}_{k\in\Nat}$ the subset of
nonzero elements of $\{a_n\}_{n\in\Nat}$. Then,

$$\gamma(T_wf)=\sum_{k=0}^{\infty}a_{n_k}\gamma(u_{n_k}^{v_{n_k}})=\sum_{k=0}^{\infty}\|\gamma_{n_k}(f)\|w^{n_k}
\gamma(u_{n_k}^{v_{n_k}})=\sum_{k=0}^{\infty}\gamma(u_{n_k}^{\gamma_{n_k}(f)})w^{n_k}=
\sum_{k=0}^{\infty}\gamma_{n_k}({\gamma_{n_k}(f)})w^{n_k}.$$

\nd where the last sum comes from definition of $\gamma$ and is
precisely $B(f,g;w)$, that is

\begin{equation}\label{eqqq1}
\gamma(T_wf)=B(f,g;w)
\end{equation}

\nd and so, for $0\leq r<1$ we obtain $|B(f,g;r)|\leq \|\gamma\|\|T_rf\|\leq \|\gamma\|A_4(H)\|f\|.$ Then $g\in H^b$ and $\Ge\leq A_4(H)\|\gamma\|$, completing the
proof. $\bl$

\begin{defn}\label{def5.8}{\rm Let $H$ be a Banach space
of type $\mathcal{H}(\mathbb{D};F)_4$. With the notation of Proposition \ref{teo12} let
$\Gamma\colon H'\longrightarrow H^b$
be defined by $\Gamma(\gamma)=g$. The mapping $\Gamma$ is clearly linear.}
\end{defn}

\begin{prop}\label{eight}{\rm Let $H$ be a Banach space of type $\holo_4$
satisfying $P_6$. Then $\Gamma\colon H^\prime  \longrightarrow H^b$ is a isometric isomorphism.}
\end{prop}

\p The proof in the scalar-valued case applies. See Taylor \cite[Theorem 9.3]{Taylor 1}. $\bl$

\hs The next result has a great importance when we study duality
of Hardy spaces:

\begin{thm}\label{thm1}{\rm Let $H$ be a Banach space of type $\holo_4$ satisfing that for all $f\in H$

$$\dss\lim_{r\rightarrow
1}\|T_rf-f\|=0.$$

Then, every $\gamma\in H'$ can be represented as follows:

\begin{equation} \label{reprf}\gamma(f)=\lim_{r\rightarrow
1}\ds{1}{2\pi}\int_{0}^{2\pi}<f(\rho
e^{i\theta}),g(\left(\frac{r}{\rho}\right)e^{-i\theta})>d\theta,
\end{equation}

\nd where $0\leq r<\rho<1$ and $g\in H^b$. The element $g\in H^b$
uniquely determines and is uniquely determined by $\gamma$ and $\|\gamma\|=\Ge$.}
\end{thm}

\p The proof in the scalar-valued case applies. See Taylor \cite[Theorem 10.1]{Taylor 1}. $\bl$

\begin{remark}{\rm From now on, if $H=H^p(\mathbb{D};F)$, we assign, for each $g\in\mathcal{H}(\mathbb{D};F')$,

$$N_{q}(g;r)=\sup_{\|f\|_p=1}|B(f,g;r)|$$

\nd instead of $N(g;z)$ and, if $g\in $\hp$^b$, we write

$$\|g\|_q^b=N_q(g)=\sup_{0\leq r<1}N_{q}(g;r).$$}
\end{remark}

\begin{prop}\label{62}{\rm Let $1\leq p\leq \infty$.
If $g\in{\cal H}(\mathbb{D};F')$, with $q$ being the conjugate indice of $p$,
we have

\begin{equation}\label{ttt4}
N_q(g;r)\leq \|g\|_{q,r}
\end{equation}

\nd which implies $H^q(\mathbb{D};F')\subset H^{p}(\mathbb{D};F)^b$.}

\end{prop}

\p The proof in the scalar-valued case applies. See Taylor \cite[Theorem 16.2]{Taylor 2}. $\bl$

\section{Boundary values and the analytic Radon-Nikodym property}

\begin{defn}{\rm Let $1\leq p\leq \infty$, and let $H^p(\mathbb{T};F)$ be defined by

\begin{center}$H^p(\mathbb{T};F)=\{{f}\in L^p(\mathbb{T};F);\displaystyle\frac{1}{2\pi}
\int_0^{2\pi}{f}(e^{i\theta})e^{-in\theta}d\theta=0, n=-1,-2,
\ldots\}.$
\end{center}

\nd Clearly $H^p(\mathbb{T}; F)$ is a closed subspace of $L^p(\mathbb{T};F)$.}
\end{defn}

\begin{defn}{\rm Let $\varphi\in L^1(\mathbb{T};F)$ and $f$ be defined by

\begin{equation}\label{ci} f(z)=\ds{1}{2\pi i}\int_\mathbb{T} \ds{\varphi(\xi)}{\xi-z}d\xi,\,|z|<1.\end{equation}

Of course $f\in H^p(\mathbb{D};F)$. We say $f$ is the Cauchy integral of $\varphi$. }
\end{defn}

The following result of Ryan extends a classical theorem of Riesz in the case of scalar-valued functions.

\begin{thm}{\rm{\bf (Ryan \cite{Ryan,Ryan 2})} \label{ryan1} Let $1 \leq p\leq \infty$, and let $F$ be a separable and reflexive Banach space. For each
$\varphi\in H^p(\mathbb{T};F)$ let $f\colon \mathbb{D}\longrightarrow F$ be defined by the Cauchy integral (\ref{ci}).

\nd Then the mapping $\varphi\longrightarrow f$ is an isometric isomorphism between $H^p(\mathbb{T};F)$ and $H^p(\mathbb{D};F)$.
Moreover the limit

\begin{center}$\varphi(e^{i\theta})=\dss\lim_{r\rightarrow
1}f(re^{i\theta})$\end{center}

\nd exists almost everywhere in the norm of $F$}
\end{thm}

\hs Ryan's result was improved by Bukhvalov as follows.

\begin{thm}\label{ARNP}
{\rm{\bf (Bukhvalov \cite{Bukhvalov})} \label{t1} Let  $1\leq p\leq
\infty$. For each $\varphi\in H^p(\mathbb{T};F)$ let $f\colon \mathbb{D}\longrightarrow F$ be defined by
the Cauchy integral (\ref{ci}). Then the mapping $\varphi\mapsto f$ is an isometric isomorphism
between $H^p(\mathbb{T};F)$ and the closed subspace of all $f\in H^p(\mathbb{D};F)$ such that the radial limits

\begin{equation}\label{5.2}\varphi(e^{i\theta})=\dss\lim_{r\rightarrow
1}f(re^{i\theta}) \end{equation}

\nd exist almost everywhere in the norm of $F$. If $\varphi_r(e^{i\theta})=f(re^{i\theta})$ and $1\leq p<\infty$, then

\begin{equation} \label{lemARNP}\varphi_r\longrightarrow \varphi \text{ in } H^p(\mathbb{T};F).\end{equation}}
\end{thm}

\begin{remark}{\rm Motivated by the last theorem, from now on if $\varphi\in H^p(\mathbb{T};F)$ and $f$ is the Cauchy integral of $\varphi$, we will
 write $\tl{f}$ instead of $\varphi$.}
\end{remark}

\begin{defn}{\rm Following Bukhvalov \cite{Bukhvalov} we say that a Banach space $F$ has the analytic Radon-Nikodym property (ARNP for short) if for each
 $1\leq p\leq \infty$, every $f\in H^p(\mathbb{D};F)$ satisfies (\ref{5.2}). Thus, if $F$ has ARNP and $1\leq p\leq \infty$, then the mapping $\tl{f}\in H^p(\mathbb{T};F)
\longrightarrow f\in H^p(\mathbb{D};F)$ is an isometric isomorphism. We say {\it $\tilde{f}$ is the boundary value function
associated with $f$.}}
\end{defn}

We recall that a Banach space $F$ has the Radon-Nikodym property (RNP for short) if for each finite measure space $(X,\Sigma,\mu)$ and each measure $\nu\colon\Sigma\longrightarrow
F$ of finite variation and absolutely continuous with respect to $\mu$, there exists $f\in L^1(X,\Sigma, \mu;F)$ such that $\nu(A)=\int_Afd\mu$ for every $A\in \Sigma$. Every
reflexive Banach space has the RNP. We refer
to the book of Diestel and Uhl \cite{Diestel} for background information on the RNP.

By Ryan's result every separable and reflexive Banach space has the ARNP. Bukhvalov \cite{Bukhvalov} has proved that every dual Banach space with the RNP has also the ARNP.
Bukhvalov and Danilevich \cite{Bukhvalov2} have proved that a Banach space $F$ has the ARNP if for some $1\leq p\leq \infty$, every $f\in H^p(\mathbb{D};F)$ satisfies (\ref{5.2}).
They also have proved that every Banach space with the RNP has also the ARNP. Burkholder \cite{Burkholder 2} has shown that $L^1(\mathbb{T})$ has the ARNP, but it does not have
the RNP.

\begin{thm}\label{ttt6}{\rm Let $F$ be a Banach space with the ARNP. Then \hp\,, for $1\leq
p<\infty$, is a Banach space of type ${\cal H}(\mathbb{D};F)_7$.}
\end{thm}

\p It is enough to prove that $P_6$ and $P_7$ are satisfied. Property $P_5$ follows then from Proposition \ref{lem41} and Proposition \ref{prop31}.

\hs Property $P_7$ follows from equation (\ref{ig2}), i.e.,
$\|T_rf\|_p=\|f\|_{p,r}$.

\hs By definition, the boundary value related, from the
isomorphism of Theorem \ref{ARNP}, to $T_rf$ is ${f}_r$, where
${f}_r(e^{i\theta})=f(re^{i\theta})$. From
(\ref{lemARNP}), we have that

$$\lim_{r\rightarrow 1}\|{f}_r-\tl{f}\|_p=0$$

\nd and again from the isomorphism of Theorem \ref{ARNP},

$$\|{f}_r-\tl{f}\|_p=\|T_r{f}-{f}\|_p$$

\nd since the boundary value of $T_rf-f$ is ${f}_r-\tl{f}$ for all
$0\leq r<1$. Therefore, $P_6$ is satisfied. $\bl$

\begin{thm}\label{teo81}{\rm If $F$ has the
ARNP and $1\leq p< \infty$, then $H^p(\mathbb{D};F)'$ is isometrically
isomorphic to $H^p(\mathbb{D};F)^b$.}
\end{thm}

\p It follows from Propositions \ref{eight} and \ref{ttt6}. $\bl$

\section{Duality and the UMD-property (UMDp)}

\hspace{1.5cm} In this section $p,q\in(1,\infty)$ denote conjugate indices. In Proposition \ref{62} we saw that $H^q(\mathbb{D};F')\subset
H^{p}(\mathbb{D};F)^b$ and we wonder when

\begin{equation}\label{81}
H^{p}(\mathbb{D};F)^b\subset H^q(\mathbb{D};F'),
\end{equation}

\nd that is, $H^{p}(\mathbb{D};F)^b= H^q(\mathbb{D};F')$.  The answer is
not always positive. In the case of a positive answer, we prove
that $H^q(\mathbb{D};F')$ and $H^p(\mathbb{D};F)'$ are topologically
isomorphic spaces by Theorem \ref{teo81}. Let us express the inclusion (\ref{81}) as an
assertion:

\vs

\nd {\bf A$_1$($q;F'$):} If $g\in H^p(\mathbb{D};F)^b$, then $g\in
H^q(\mathbb{D};F')$.

\begin{prop}\label{82}{\rm Let $1<q<\infty$. If  {\bf A$_1$($q;F'$)} is true, then there
exists a constant $C_1(q)>0$, depending only on $q$, satisfying the
inequality

\begin{equation}\label{qq1}\|g\|_{q,r}\leq C_1(q)N_q(g;r),
\end{equation}

\nd for every $0\leq r<1$ and $g\in {\cal H}(\mathbb{D};F')$. In
particular, the inclusion given in (\ref{81}) is continuous.}

\end{prop}

\p The proof in the scalar-valued case applies. See Taylor \cite[Theorem 17.1]{Taylor 2}. $\bl$

\vs

Consider the assertion:

\nd {\bf A$_2$($p;F$):} Let $f\in \mathcal{H}(\mathbb{D};F)$ be the
Cauchy integral of $\varphi\in L^p(\mathbb{T};F)$. Then $f\in
H^p(\mathbb{D};F)$.

\begin{prop}\label{pff3}{\rm If {\bf A$_2$($p;F$)} is true, then there exists a constant $C_2(p)>0$, depending only on $p$, such that:

$$\|f\|_p\leq C_2(p) \|\varphi\|_p,$$

\nd where $\varphi\in L^p(\mathbb{T};F)$ and $f\in H^p(\mathbb{D};F)$ is the
Cauchy integral of $\varphi$.

}
\end{prop}

\p The proof in the scalar-valued case applies. See Taylor \cite[Theorem 17.2]{Taylor 2}. $\bl$

\vs

\hs We will now establish some implications between the assertions ${\bf A_1}(q; F')$ and ${\bf A_2} (p; F)$. These assertions will
be useful in Theorem \ref{teo 111} when we characterize the dual of $H^p(\mathbb{D};F)$ under certain conditions on $F$.

\begin{prop}\label{prop b}{\rm Let $\tl{g}\in L^{q}(T;F')$, $1\leq
p\leq \infty$ and $0\leq r<1$. If $g$ is the Cauchy Integral of $\tl{g}$, then
$g\in H^p(\mathbb{D};F)^b$ and $\|g\|_q^b\leq \|\tl{g}\|_{q}$.}
\end{prop}

\p From the
fact that the Bochner integral commutes with linear functionals
(see Mujica \cite[page 42]{Mujica}), the proof given in the scalar-valued case applies. See Taylor \cite[page 32, Theorem 14.1]{Taylor 2}.
$\bl$

\begin{prop}\label{pf5}{\rm {\bf A$_1$($q;F'$)} implies {\bf A$_2$($q;F'$)}.}
\end{prop}

\p Let $\varphi \in L^q(T;F')$ and $g$ be the Cauchy integral of
$\varphi$. By Proposition \ref{prop b}, $g\in H^p(\mathbb{D};F)^b$.
Since {\bf A$_1$($q;F'$)} is true, $g\in H^q(\mathbb{D};F')$. $\bl$

\begin{prop}\label{pff2}{\rm {\bf A$_2$($p;F$)} implies {\bf A$_1$($q;F'$)}.}
\end{prop}

\p Let $g\in H^p(\mathbb{D};F)^b$, $\varphi\in L^p(\mathbb{T};F)$ and $f$ be
the the Cauchy integral of $\varphi$. In this way, since $f$ is
the Cauchy integral of $\varphi$:

$$\ds{1}{2\pi}\int_0^{2\pi}<\varphi(e^{i\theta}),g(re^{-i\theta})>d\theta=\sum_{n=0}^{\infty}
<c_n(\varphi),\gamma_n(g)>r^n=B(f,g;r),$$

\nd where
$c_n(\varphi)=\frac{1}{2\pi}\int_0^{2\pi}\varphi(e^{i\theta})e^{-in\theta}d\theta$,
for all $n\in \Int$. By Proposition \ref{pff3}:

$$\left|\ds{1}{2\pi}\int_0^{2\pi}<\varphi(e^{i\theta}),g(re^{-i\theta})>d\theta\right|\leq \|f\|_p\|g\|_q^b\leq
C_2(p)\|\varphi\|_p\|g\|_q^b,$$

\hs By a result of Diestel and Uhl \cite[page 97]{Diestel},

$$\|g\|_{q,r}\leq C_2(p)\|g\|_q^b$$

\nd and the proof is complete. $\bl$

The unconditional martingale difference property (UMDP for short) has been extensively studied by Burkholder \cite{Burkholder 2}.
We say that $F$ has the UMDP or that $F$ is UMD.

Instead of giving the original definition of UMDP, for which it would be necessary to introduce several other concepts, we will give a characterization.
First, we need the next definition:

\begin{defn}{\rm If $f\in L^p(\mathbb{T};F)$, the {\it analytic
projection $f^a$ of $f$} is the function whose negative Fourier
coefficients are zero and the other coincide with the respective
Fourier coefficients of $f$.}\end{defn}

It is known the Hilbert transform with values in $F$ is bounded in $L^p(T;F)$ if and only if the analytic projection is bounded. This can be proved using
Hoffman \cite[page 151]{Hoffman}, Riesz \cite{Riesz} and Burkholder \cite[pages 408 and 409]{Burkholder 3}. It is also known that the boundedness of Hilbert
transform as an operator on $L^p(T;F)$ is equivalent to $F$ being UMD. The sufficiency was established
by Burkholder \cite{Burkholder 3, Burkholder} and the necessity by Bourgain \cite{Bourgain}. So, we have the following theorem:

\begin{thm}\label{umdp2}\,{\rm A Banach space $F$ is UMD
if and only if

$$\begin{array}{rcrcclcccc}
   &S\colon L^p(\mathbb{T};F)& \longrightarrow & H^p(\mathbb{T};F) \\
  &f& \longmapsto & f^a
\end{array}$$

\nd is a bounded linear map for every $1<p<\infty$, where $f^a$ is
the analytic projection of $f$.}
\end{thm}

\begin{remark}\label{rem6.5} {\rm An interesting result is that any
UMD-space is reflexive, in fact supperreflexive (see Aldous \cite{aldous}
or Maurey \cite{Maurey}). In particular, if $F$ is UMD, then $F$ and $F'$ have the RNP. On the other hand Pisier \cite{Pisier} has constructed an example showing that a supperreflexive space
need not be UMD.}
\end{remark}

\begin{prop} \label{A2}{\rm A Banach space $F$ with the ARNP is UMD if and only
if it satisfies ${\bf A}_2(p;F)$.}
\end{prop}

\p In this proof we use the characterization of UMD-property given
by Theorem \ref{umdp2}. Suppose {\bf A$_2$($p;F$)} is true. Then if
$\varphi\in L^p(\mathbb{T};F)$, the Cauchy integral $f$ of $\varphi$ is an
element of $H^p(\mathbb{D};F)$ and by Bukhvalvov
\cite[Theorem 2.3]{Bukhvalov}, ${\varphi}$ is the unique element of
$L^p(\mathbb{T};F)$ such that $f$ is the Cauchy integral. By Theorem
\ref{ARNP}, there is only one $\tl{f}\in H^p(\mathbb{T};F)\subset L^p(\mathbb{T};F)$
where $f$ is the Cauchy integral of $\tl{f}$. It is not difficult
to see that $\tl{f}$ is the analytic projection of $\varphi$. Thus, $F$ is UMD by Proposition \ref{pff3}.

\hs Now suppose $F$ is UMD, that is, by Theorem \ref{umdp2} the linear map:

$$\begin{array}{rcrcclcccc}
   &S\colon L^p(\mathbb{T};F)& \longrightarrow & H^p(\mathbb{T};F) \\
  &\varphi& \longmapsto & \varphi^a
\end{array},$$

\nd is bounded for all $1<p<\infty$, i.e., there is $C>0$ such that $\|\varphi^a\|_p\leq C\|\varphi\|_p$. If $f$ is the Cauchy
integral of $\varphi$, for $z\in \mathbb{D}$ we have:

$$f(z)=\ds{1}{2\pi}\int_0^{2\pi}\ds{\varphi(e^{i\theta})}{1-ze^{i\theta}}d\theta=
\ds{1}{2\pi}\int_0^{2\pi}\ds{\varphi^a(e^{i\theta})}{1-ze^{i\theta}}d\theta$$

\nd and by Theorem \ref{ARNP}, $\|f\|_p=\|\varphi^a\|_p$. This concludes the proof. $\bl$

\hyphenation{de-mons-tra-ção}

\begin{prop}\label{ppp1}{\rm Let $F$ be a Banach space that has the ARNP. Then, {\bf
A$_1$($q;F'$)} implies that $F$ is UMD. Furthermore, if $F'$ also has the ARNP, the converse implication is
true.}
\end{prop}

\p If $F$ is UMD, then {\bf A$_1$($q;F'$)} is true by Proposition \ref{pff2} and Proposition
\ref{A2}. Now, if {\bf A$_1$($q;F'$)} is
true, by Proposition \ref{pf5} {\bf
A$_2$($q;F'$)} is true, and so $F'$ is UMD by Proposition \ref{A2}. By
results of Rubio de Francia \cite[page 205]{Rubio} and Burkholder
\cite[page 237]{Burkholder 2}, we obtain that $F$ is UMD. $\bl$

\begin{prop}\label{teof1}{\rm  If $F$ has the ARNP and {\bf A$_1$($q;F'$)}
is true, then the spaces $H^p(\mathbb{D};F)'$ and $H^q(\mathbb{D};F')$ are
topologically isomorphic.}
\end{prop}

\p If $g\in H^q(\mathbb{D};F')$, then from the inequalities (\ref{ttt4})
and (\ref{qq1}), it follows that:

$$\|g\|_q^b\leq \|g\|_q\leq C_1(q)\|g\|_q^b.$$

\hs The conclusion follows by Theorem \ref{teo81}. $\bl$

\begin{defn}{\rm We say that $H^p(\mathbb{D};F)'$ and $H^q(\mathbb{D};F')$
are canonically topologically isomorphic if {\bf A$_1$($q;F'$)} is
true.}
\end{defn}

\begin{thm}\label{teo 111}{\rm Let $F$ be a Banach space with ARNP such that $F'$ also has the ARNP.
Then $H^p(\mathbb{D};F)^\prime$ and $H^q(\mathbb{D};F')$ are canonically topologically
isomorphic if and only if $F$ is UMD. In this case, we have the isomorphism

$$\begin{array}{rccrcll}\Psi_p\colon & H^q(\mathbb{D};F')&\longrightarrow &
H^p(\mathbb{D};F)'\\  &g&\longmapsto & \Psi_p(g)\colon&
H^p(\mathbb{D};F)&\longrightarrow &\Complex\\ & & & &f& \longmapsto
&\Psi_p(g)(f)=\ds{1}{2\pi}\int_0^{2\pi}<\tl{f}(e^{i\theta}),\tl{g}(e^{-i\theta})>d\theta,
\end{array}$$

 \nd where $\tl{f}$ and $\tl{g}$ are the boundary values of $f\in H^p(\mathbb{D};F)$
 and $g\in H^q(\mathbb{D};F)$, respectively.}
\end{thm}

\p The theorem follows from Propositions \ref{ppp1} and \ref{teof1} with the aid of Theorems \ref{thm1} and \ref{ttt6}. $\bl$

\begin{remark}{\rm For similar results see Bukhvalov \cite[Theorem 3.1]{Bukhvalov} and Blasco \cite{Blasco}.}
 \end{remark}

\begin{cor}\label{reflexive}{\rm If $F$ is UMD, then $H^p(\mathbb{D};F)$
is reflexive when $1<p<\infty$.}
\end{cor}

\p We consider the embedding $J\colon H^p(\mathbb{D};F)\longrightarrow H^p(\mathbb{D};F)''$ and, since $F$ is reflexive (see Remark
\ref{rem6.5}), the isomorphisms we have obtained
in Theorem \ref{teo 111}:

\begin{center}$\Psi_p\colon
H^q(\mathbb{D};F')\longrightarrow H^p(\mathbb{D};F)'$ and $\Psi_q\colon
H^p(\mathbb{D};F)\longrightarrow H^q(\mathbb{D};F')'$.\end{center}

\hs Given $x''\in H^p(\mathbb{D};F)''$, we have $x''\circ\Psi_p\in H^q(\mathbb{D};F')'$. Thus, there is $f\in H^p(\mathbb{D};F)$ such that
$\Psi_q(f)=x''\circ \Psi_p$. Also, if $x'\in H^p(\mathbb{D};F)'$ there is $g\in H^q(\mathbb{D};F')$ such that $\Psi_p(g)=x'$. Then, it follows that

$$x''(x')=x''(\Psi_p(g))=\Psi_q(f)(g)=\frac{1}{2\pi}\int_0^{2\pi}\hspace{-.2cm}<\tl{f}(e^{i\theta}),\tl{g}(e^{-i\theta})>d\theta=
\Psi_p(g)(f)=x'(f)=J(f)(x')$$

\nd and the proof is complete. $\bl$

\section{Weak spaces of vector-valued functions}

If $A$ is a subset of some Banach space
$E$ and $\mathcal{F}(A;F)$ is a set of functions $f\colon
A\longrightarrow F$ satisfying some property, then the weak version of $\mathcal{F}(A;F)$
is given by

$$\mathcal{F}(A;F)_w=\{f\colon A\longrightarrow F; \psi\circ f\in
\mathcal{F}(A;\Complex), \forall \psi \in F^\prime\}.$$

As in the case of Hardy and Lebesgue spaces, when
$F=\Complex$ we write $\mathcal{F}(A)$ instead of
$\mathcal{F}(A;\Complex)$.

In this section we study the weak versions of the spaces $L^p(\mathbb{T};F)$
and $H^p(\mathbb{D};F)$: $L^p(\mathbb{T};F)_w$ and $H^p(\mathbb{D};F)_w$, respectively. We will see some properties of these spaces.
In the case of $H^p(\mathbb{D};F)$, we relate the study with duality. The inspiration for this section came from a result about
holomorphic functions that we can see in Mujica \cite[page 65]{Mujica}: for all Banach space $F$, $\mathcal{H}(\mathbb{D};F)=\mathcal{H}(\mathbb{D};F)_w$.

\begin{prop}\label{fjb1}{\rm For each $1\leq p<\infty$ and $f\in L^p(\mathbb{T};F)_w$ the following
supremum is finite:

\begin{center}$\dss\sup_{\psi\in F', \|\psi\|\leq
1}\|\psi\circ f\|_p$.\end{center}}
\end{prop}

\p We prove the case $1\leq p< \infty$. For
$p=\infty$, the proof is analogous. For each $f\in
L^p(\mathbb{T};F)_w$, consider the linear map given by:

$$\begin{array}{rll} S_{f}\colon F'& \longrightarrow & L^p(\mathbb{T}) \\
\psi& \longmapsto& \psi\circ f
\end{array}$$

\hs We prove that $S_{f}$ is continuous using the Closed Graph Theorem. From this, the result follows. In fact:

\begin{center} $\|S_{f}\|=\dss\sup_{\psi\in F', \|\psi\|\leq
1}\left(\ds{1}{2\pi}\int_0^{2\pi}|\psi\circ
f(e^{i\theta})|^pd\theta\right)^p$.
\end{center}

\hs Let $(\psi_n)_{n\in\Nat}\subset F'$ be such that $\psi_n\longrightarrow \psi\in F'$.
Now suppose $S_{f}(\psi_n)\longrightarrow g\in L^p(\mathbb{T})$, that is,

\begin{center}$\dss\lim_{n\rightarrow
\infty}\left(\ds{1}{2\pi}\int_0^{2\pi}\|\psi_n\circ f(e^{i\theta})
-g(e^{i\theta})\|^pd\theta\right)^{\frac{1}{p}}=0$.\end{center}

\hs By Bartle \cite[page 73]{Bartle}, there exists a subsequence $(\psi_{n_k}\circ f)$
such that ${\psi_{n_k}\circ f}(e^{i\theta})\longrightarrow g(e^{i\theta})$ almost everywhere. Thus, also almost everywhere

$$g(e^{i\theta})=\psi\circ f(e^{i\theta}),$$

\nd that is, $g=S_{f}(\psi)\in L^p(\mathbb{T})$. $\bl$

\begin{defn}{\rm For each $f\in L^p(\mathbb{T};F)_w$ and $1\leq p\leq
\infty$, we write

\begin{center}$\|f\|^w_p=\dss\sup_{\psi\in F', \|\psi\|\leq
1}\|\psi\circ f\|_p$.\end{center}}
\end{defn}

\begin{prop}\label{1.13}{\rm Let $1\leq p\leq \infty$, let $F$ be a Banach space and suppose that there exists a sequence $(\psi_n)\subset F^\prime$ which separates
the points of $F$. Then:

 \nd {\bf a)} $(L^p(\mathbb{T};F)_w, \|.\|_p^w)$ is a normed space.

\nd {\bf b)} $(L^p(\mathbb{T};F)_w,\|.\|_p^w)$ is isometrically isomorphic to a subspace of ${\cal L}(F';L^p(\mathbb{T}))$.

\nd {\bf c)} $(L^p(\mathbb{T};F)_w, \|.\|_p^w)$ is a Banach space.}
\end{prop}

\p It is easy to see that
$(L^p(\mathbb{T};F); \|.\|^w_p)$ is a seminormed space for every $F$. If $f\in L^p(\mathbb{T};F)$ and $\|f\|^w_p=0$, then it follows that $\psi_n\circ f=0$ a.e. for every
 $n\in \Nat$. Since $(\psi_n)$ separates the points of $F$, it follows that $f=0$ a.e.

For assertion (b), define the map

$$\begin{array}{ccccc} S\colon L^p(\mathbb{T};F)_w&\longrightarrow &{\cal L}(F';L^p(\mathbb{T})) \\
  f & \longmapsto & S_f,\end{array}$$

\nd where $S_f$ was defined in the proof of Proposition \ref{fjb1}. Of course $S$ is well-defined and is a isometric map:

$$\|S_f\|=\sup_{\psi\in F', \|\psi\|\leq
1}\|S_f(\psi)\|=\sup_{\psi\in F', \|\psi\|\leq 1}\|\psi\circ
f\|_p=\|f\|_p^w.$$

\hs Also, $S$ is injective and linear.
Assertion (c) is a consequence of assertion (b) and of the fact that $\mathcal{L}(F';L^p(\mathbb{T}))$ is a Banach space.

\begin{remark}{\rm There exists a sequence $(\psi_n)\subset F'$ which separates the points of $F$ if and only if exists an injective operator $T\in
 L(F;l_\infty)$. In particular every separable Banach space verifies this condition.}
\end{remark}

\begin{thm} {\rm For every Banach space $F$ the following conditions are equivalent.

\nd {\bf a)} $L^p(\mathbb{T};F)=L^p(\mathbb{T};F)_w$ for all $1\leq p<\infty$.

\nd {\bf b)} $L^p(\mathbb{T};F)=L^p(\mathbb{T};F)_w$ for some $1\leq p<\infty$.

\nd {\bf c)} $F$ has finite dimension.}

\end{thm}

\p Of course  (a)$\Longrightarrow $ (b). For (b)$\Longrightarrow $ (c), consider a sequence
$(x_n)_{n\in\Nat}\subset F$ such that

$$\dss\sum_{n=1}^{\infty}|\varphi(x_n)|^p<\infty$$

\nd for all $\varphi\in F'$ and $1\leq p<\infty$. Our aim is to show that
$\dss\sum_{n=1}^\infty\|x_n\|^p<\infty$, i.e., to prove that the identity $I\colon F\longrightarrow F$ is absolutely $p$-summing. From this,
by Pietsch \cite{Pietsch}, $F$ has finite-dimension.

\hs Define

\begin{center}${f}(e^{i\theta})=\dss\sum_{n=1}^{\infty}\ds{x_n}{2^{\frac{-(n+1)}{p}}\pi^{\frac{-1}{p}}}\chi_{I_n}(e^{i\theta}),$\end{center}

\nd where $0\leq \theta \leq 2\pi$, $\{I_n\}_{n=1}^{\infty}$ is a sequence of disjoint intervals, $\chi_{I_n}$ is the
characteristic function of $I_n$ and $\frac{1}{2\pi}\int_{I_n}d\theta=\frac{1}{2^{n}}$. We assert that ${f}\in L^p(\mathbb{T};F)_w$. Indeed, it is
plain that for each $\varphi \in F'$, $\varphi\circ \tl{f}$ is measurable and

$$\left(\ds{1}{2\pi}\int_0^{2\pi}|\varphi({f}(e^{i\theta}))|^pd\theta
\right)^{\frac{1}{p}}=\left(\sum_{n=1}^{\infty}\ds{1}{2\pi}\int_{I_n}|\varphi({f}(e^{i\theta}))|^p
d\theta\right)=\left(\dss\sum_{n=1}^{\infty} |\varphi(x_n)|^p
\right)^{\frac{1}{p}}<\infty.$$

\hs By hypothesis, ${f}\in L^p(\mathbb{T};F)$ and

$$\left(\sum_{n=1}^{\infty}\|x_n\|^p\right)^{\frac{1}{p}}=\left(\sum_{n=1}^{\infty}\ds{1}{2\pi}\int_{I_n}
\|{f}(e^{i\theta})\|^pd\theta\right)^{\frac{1}{p}}=
\left(\ds{1}{2\pi}\int_0^{2\pi}\|{f}(e^{i\theta})\|^pd\theta\right)^{\frac{1}{p}}
<\infty.$$

\nd Now, let us prove that (c)$\Longrightarrow $ (a). We just have to show
$L^p(\mathbb{T};F)_w\subset L^p(\mathbb{T};F)$. The other inclusion is obvious for all $F$. Let
$\{e_1,\ldots, e_n\}$ be a basis of $F$ with dual basis
$\{\varphi_1,\ldots, \varphi_n\}$.
Let $f\in L^p(\mathbb{T};F)_w$ and write $f=f_1e_1+\ldots+f_ne_n$, where $f_j=\varphi_j\circ f$ for $j=1,\ldots,n$.
Of course each $f_j\in L^p(\mathbb{T})$. Furthermore,

$$\left(\ds{1}{2\pi}\int_0^{2\pi}|f_j(e^{i\theta})|^pd\theta\right)^\frac{1}{p}
=\left(\ds{1}{2\pi}\int_0^{2\pi}|\varphi_j(f(e^{i\theta}))|^pd\theta)\right)^{\frac{1}{p}}<\infty$$

\nd resulting that $f\in L^p(\mathbb{T};F)$. $\bl$

The proof of this theorem follows a suggestion of Professor Oscar Blasco.

\begin{prop}\label{fjb2}{\rm If $1\leq p\leq \infty$ and $f\in H^p(\mathbb{D};F)_w$, the following supremum is finite:

$$\sup_{\psi\in F', \|\psi\|\leq 1}\|\psi\circ f\|_p.$$}

\end{prop}

\p If $f\in H^p(\mathbb{D};F)_w$, define

$$\begin{array}{rll} \Xi_f\colon F'& \longrightarrow & H^p(\mathbb{D}) \\
\psi& \longmapsto& \psi\circ f
\end{array}.$$

\hs Let $(\psi_n)_{n\in\Nat}\subset F'$ such that $\psi_n\longrightarrow\psi\in F'$. Suppose $\Xi_f(\psi_n)\longrightarrow g\in
H^p(\mathbb{D})$ i.e.,

$$\lim_{n\rightarrow \infty}\sup_{0\leq r<1}\left(\ds{1}{2\pi}\int_0^{2\pi}
\|\psi_n\circ
f(re^{i\theta})-g(re^{i\theta})\|^pd\theta\right)^{\frac{1}{p}}=0.$$

\hs So, for each $r$ we have

$$\lim_{n\rightarrow \infty}\left(\ds{1}{2\pi}\int_0^{2\pi}
\|\psi_n\circ
f(re^{i\theta})-g(re^{i\theta})\|^pd\theta\right)^{\frac{1}{p}}=0$$

\nd and proceeding as in the proof of Proposition
\ref{fjb1}, we obtain almost everywhere

\begin{equation}\label{fjb3}g(re^{i\theta})=\psi \circ
f(re^{i\theta}).
\end{equation}

\hs Since $g$, $f$ e $\psi$ are continuous (in particular, $f\in \mathcal{H}(\mathbb{D};F)$), for each $0\leq r< 1$, $g(re^{i\theta})=\psi \circ
f(re^{i\theta})$ for all $0\leq \theta\leq 2\pi$. Thus $g=\psi\circ f$. By the Closed Graph Theorem, $\Xi_f$
is continuous and the result follows. For $p=\infty$, the proof is analogous and we can use Lemma \ref{nonincreasinglemma}.$\bl$

\begin{defn}{\rm For each $f\in H^p(\mathbb{D};F)_w$ and $1\leq p\leq
\infty$, we write

$$\|f\|_p^w=\sup_{\psi\in F', \|\psi\|\leq 1}\|\psi\circ f\|_p.$$}\end{defn}

\begin{prop}\label{hardy215}{\rm Let $1\leq p\leq \infty$. Then:

\nd {\bf a)}($H^p(\mathbb{D};F)_w$, $\|.\|_p^w$) is a Banach space of type ${\cal H}(\mathbb{D};F)_4$ and
$(H^\infty(\mathbb{D};F)_w,\|.\|_\infty^w)=(H^\infty(\mathbb{D};F),\|.\|_\infty)$.

\nd {\bf b)} If $F$ is reflexive, then $(H^p(\mathbb{D};F)_w,\|.\|_p^w)$ is isometrically isomorphic to ${\cal L}(F';H^p(\mathbb{D}))$.}

\end{prop}

\p If $f\in H^p(\mathbb{D};F)_w$ and $\|f\|_w^p=0$, note that $\psi\circ f=0$ a.e. for all $\psi \in F^\prime$. Since $f$ and $\psi$ are continous, then
$\psi\circ f\equiv 0$, for all $\psi\in F^\prime$, which implies $f\equiv 0.$ The other properties of norm are easy to prove.

\hs To show the fact that is a Banach space of type $\mathcal{H}(\mathbb{D};F)_4$ we use a consequence of the Hahn-Banach Theorem, that asserts
for every $v\in F$ that

$$\|v\|=\sup_{\|\psi\|\leq 1, \psi\in F'}|\psi(v)|,$$

\nd and the steps of the proof of Taylor \cite[Theorem 11.1]{Taylor 2}. See also Proposition \ref{prop31}. The second part of (a) is a consequence of the Banach-Steinhaus Theorem.

For assertion
(b), define the map

$$\begin{array}{ccccc} \Xi\colon H^p(T;F)_w&\longrightarrow &{\cal L}(F';H^p(\mathbb{D})) \\
  f & \longmapsto & \Xi_f,\end{array}$$

\nd where $\Xi_f$ was defined in Proposition \ref{fjb2}. Of course $\Xi$ is well-defined and is an isometric map:

$$\|\Xi_f\|=\sup_{\psi\in F', \|\psi\|\leq
1}\|\Xi_f(\psi)\|=\sup_{\psi\in F', \|\psi\|\leq 1}\|\psi\circ
f\|_p=\|f\|_p^w.$$

Also, $\Xi$ is injective and linear. We just need to show $\Xi$ is surjective. Let $R\in {\cal L}(F';H^p(\mathbb{D}))$. Then,
for each $0\leq r<1$ and $e^{i\theta}\in T$,
we define the following element of $F''$:

$$\begin{array}{cclcl}
  \hspace{.5cm} R_{\theta,r}\colon F'&\longrightarrow & \Complex\\
  \hspace{1.4cm}\psi& \longmapsto &
  R(\psi)(re^{i\theta})
\end{array}$$

Indeed, $R_{\theta,r}$ is linear and also continuous, since is the composition of the mapping $R$, which is continuous, and of the
evaluation at the point
$re^{i\theta}$, which is continuous by Proposition \ref{teo1}.

Since $F$ is reflexive, there exists a unique $a_{\theta,r}\in F$ such that
$R_{\theta,r}(\psi)=\psi(a_{\theta,r})$, for all $\psi\in F'$. Then is well-defined the function

$$\begin{array}{cclcl}
  \hspace{.5cm} f\colon\mathbb{D} &\longrightarrow & F\\
  \hspace{1.4cm}re^{i\theta}& \longmapsto &
  a_{\theta,r}
\end{array}$$

\nd and for all $0\leq r\leq 1$, $e^{i\theta}\in \mathbb{T}$ and $\psi\in F'$ it follows that

$$\psi(f(re^{i\theta}))=\psi(a_{\theta,r})=R_{\theta,r}(\psi)=R(\psi)(re^{i\theta})$$

\nd which implies that $\psi\circ f=R(\psi)$, resulting that $f\in H^p(\mathbb{D};F)_w$ and $\Xi_f=R$.
$\bl$

\begin{thm}\label{teo fjb1}{\rm Let $p,q\in (1,\infty)$ be conjugate indices, and let
$F$ be a Banach space such that $F$ and $F'$ have the ARNP. If $H^p(\mathbb{D};F)_w=H^p(\mathbb{D};F)$, then:

\nd {\bf a)} $H^p(\mathbb{D};F)'$ and $H^q(\mathbb{D};F')$ are canonically topologically isomorphic.

\nd {\bf b)} $F$ is UMD.}
\end{thm}

\p Let us show that the equality $H^p(\mathbb{D};F)_w=H^p(\mathbb{D};F)$ implies {\bf
A$_2$($p;F$)}. The rest of the proof follows from
Propositions \ref{pff2} and \ref{teof1} and Theorem \ref{teo 111}. Let $\varphi \in
L^p(\mathbb{T};F)$ and let $f\in {\cal H}(\mathbb{D};F)$ be its Cauchy integral. For all $\psi\in F'$, we have $\psi\circ \varphi\in L^p(\mathbb{T})$ and by
Taylor \cite[page 45]{Taylor 2} we obtain $\psi\circ f\in H^p(\mathbb{D})$. By hypothesis, $f\in H^p(\mathbb{D};F)$. It concludes the proof. $\bl$

\begin{exe}\label{example}{\rm
Let $p,q\in (1,\infty)$ be conjugate indices. If $F=H^q(\mathbb{D})$, then $F$ is UMD, but $H^p(\mathbb{D};F)_w\neq H^p(\mathbb{D};F)$.}\end{exe}

\hs Suppose that $H^p(\mathbb{D};F)_w=H^p(\mathbb{D};F)$ for all
$1<p<\infty$. We know that
$H^q(\mathbb{D})$ and $H^q(\mathbb{T})$ are isometrically isomorphic. By
Pelczynski \cite[page 8]{Pelczynski},  $H^q(\mathbb{T})$ and $L^q(\mathbb{T})$ are isomorphic. By Burkholder \cite[page 237]{Burkholder 2}, every Banach
space isomorphic to a UMD-space is also UMD and $L^q(\mathbb{T})$ is UMD. Hence, we obtain that
$H^q(\mathbb{D})$ is UMD.
It follows by Proposition \ref{hardy215}, Corollary
\ref{reflexive}, and the assumption $H^p(\mathbb{D};F)_w=H^p(\mathbb{D};F)$ that ${\cal L}(F';H^p(\mathbb{D}))$ is reflexive. Let ${\cal L}_k(F';H^p(\mathbb{D}))$
denote the subspace of all
$T\in\mathcal{L}(F';H^p(\mathbb{D}))$ which are compact. Also, since $L^q(\mathbb{T})$ has the approximation property (see Diestel and Uhl
\cite[page 245]{Diestel}), it follows
that $H^q(\mathbb{D})$ has the approximation property and by Mujica
\cite[Theorem 2.1]{Mujica 2} we have that ${\cal L}(F';H^p(\mathbb{D}))={\cal L}_k(F';H^p(\mathbb{D})).$ So, if we consider $F=H^q(\mathbb{D})$ it follows that
$F'=H^p(\mathbb{D})$ and then

\begin{center}${\cal L}(H^p(\mathbb{D});H^p(\mathbb{D}))={\cal L}_k(H^p(\mathbb{D});H^p(\mathbb{D})).$\end{center}

\hs This implies that the identity $I\colon H^p(\mathbb{D})\longrightarrow
H^p(\mathbb{D})$ is compact, which is false, since
$H^p(\mathbb{D})$ is a infinite-dimensional space.  $\bl$

\hs An example where the equality $H^p(\mathbb{D};F)_w=H^p(\mathbb{D};F)$ occurs is when $F$ is finite dimensional. We have just seen an example
where the equality $H^p(\mathbb{D};F)_w=H^p(\mathbb{D};F)$ is not true, where $F$ is an infinite dimensional space.

\begin{conjecture}{\rm We conjecture that $H^p(\mathbb{D};F)=H^p(\mathbb{D};F)_w$ if and only if $F$ is finite dimensional.}\end{conjecture}

\end{document}